\documentclass[11pt]{article}

\usepackage{amssymb}
\usepackage{amsmath}
\usepackage{amsthm}
\usepackage{graphicx}
\usepackage{accents}

\newcommand{\R}{\mathbb R}
\newcommand{\C}{\mathbb C}
\newcommand{\tD}{\mathbb D}

\newcommand{\Oh}{{\cal O}}

\usepackage{xcolor,todonotes}

\author{Jens Lang}
\title{Rosenbrock-Wanner Methods: Construction and Mission}
\author{
Jens Lang\\[0.5cm]
{\small \it Technical University of Darmstadt} \\
{\small \it Department of Mathematics} \\
{\small \it Dolivostra{\ss}e 15, 64293 Darmstadt, Germany}\\
{\small lang@mathematik.tu-darmstadt.de}
}
\begin{document}
\maketitle

\begin{abstract}
This paper is concerned with the history of Rosenbrock-Wanner methods
first initiated by Rosenbrock in 1963. His original
ideas are highlighted and the main developments over 56 years
are reviewed.
\end{abstract}



\section{Introduction}
Howard Harry Rosenbrock (1920-2010) suggested in his famous paper from 1963
\cite{Rosenbrock1963} to replace the iterative process for the solution of
nonlinear problems within an implicit time integrator by a finite number of
solutions of linear systems. He summarized: {\it Some general implicit processes
are given for the solution of simultaneous first-order differential equations. These processes, which use successive substitution, are implicit analogues of the (explicit) Runge-Kutta processes. They require the solution in each time step of one or more
sets of simultaneous linear equations, usually of a special and simple form. Processes of any required order can be devised, and they can be made to have a wide margin of stability when applied to a linear problem.} Thus, Rosenbrock methods avoid the problem of
convergence for the solution of systems of nonlinear equations, making them a good alternative to fully implicit Runge-Kutta methods.
In this note, I will give a brief historical overview and explain main construction principles including widely used members of the whole family of linearly implicit methods such as Rosenbrock-Wanner methods, W-methods, and recently developed two-step Rosenbrock-Peer and W-methods. Rosenbrock methods have made their way into real-life applications and become part of very sufficient adaptive multilevel PDE-solvers, see e.g. \cite{Lang2000}. Nowadays, there still is an increasing interest in these methods, which would have delighted Rosenbrock who concluded his paper by expressing his wish: {\it The processes described above have been explored only cursorily, and it is hope that this note may stimulate others to investigate their possibilities.} It certainly did.\\

\begin{minipage}[h]{8.5cm}
Who was Howard H. Rosenbrock? Rosenbrock was born on December 16, 1920 in Ilford, England.
He graduated 1941 from University College London with a 1st class honors degree in Electrical Engineering and received his PhD from London University in 1955. During the 1960s he worked at the Cambridge University and the MIT. In 1966, he became the Chair of Control Engineering at the University of Manchester, Institute of Science and Technology.
He died on 21 October 2010.
Rosenbrock produced over $120$ scientific papers, $7$ books, and about $30$ papers on the philosophical basis of science and technology. An obituary was published in \cite{Wellstead2011}.
\end{minipage}\hspace{0.5cm}
\begin{minipage}[h]{5cm}
\includegraphics[height=5cm]{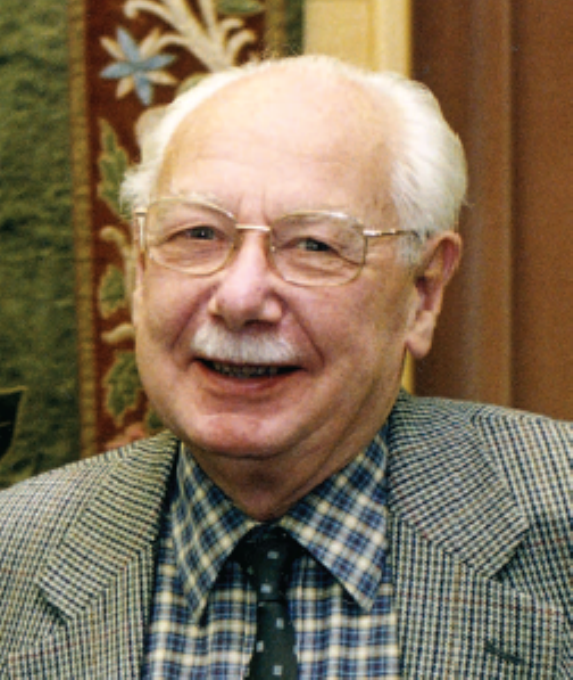}\\[2mm]
{\small\textcopyright\hspace{1mm}IEEE Control System \cite{Wellstead2011}.}
\end{minipage}

\section{The Original Idea of Rosenbrock}
In what follows, I will first review the original ideas
of Rosenbrock as described in \cite{Rosenbrock1963}. The
writing is presented in a modern style, but a few text
passages are included as pictures.

\begin{figure}[ht]
\centering
\fbox{\includegraphics[width=0.85\textwidth]{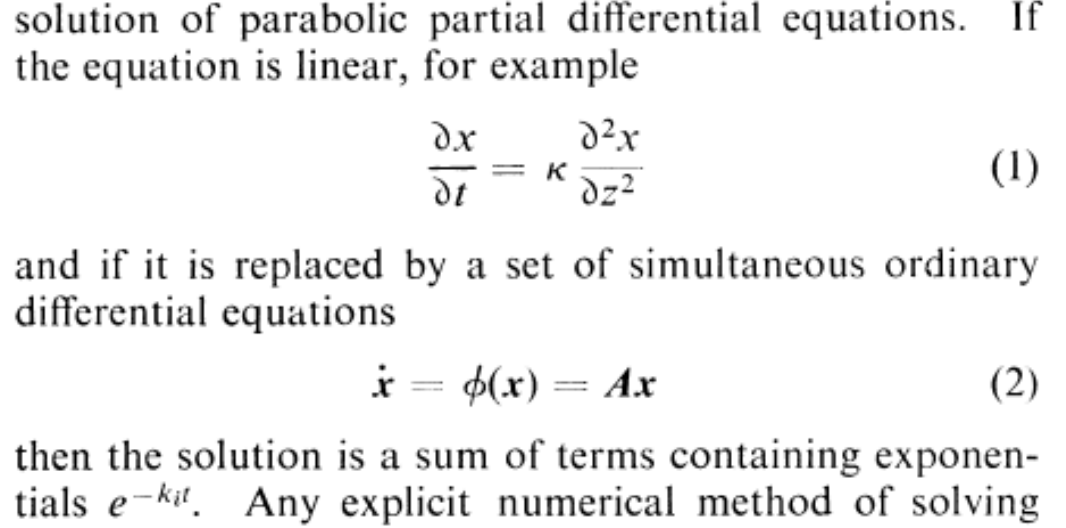}}
\parbox{13cm}{
\caption{\textcopyright\hspace{1mm}The Computer Journal, part of page 329 of \cite{Rosenbrock1963}.}
\label{fig:rosenbrock-paper-pic1}
}
\end{figure}

As starting point in his paper, Rosenbrock took a look at
the (spatial) semi-discretization of the one-dimensional linear heat
equation, i.e., formulas (1) and (2) in Fig.~\ref{fig:rosenbrock-paper-pic1}.
He stated: {\it Any explicit numerical method of solving eqn. (2) (e.g. Runge-Kutta)
replaces the exponentials by their truncated Taylor's series during one
time interval of the solution. The exponentials tend to zero as $t$
becomes large, whereas the truncated Taylor's series tend to infinity. A
severe limitation on the length of the time intervals is thus introduced.}

To illustrate these fundamental observations, let us consider the heat
equation in the form
\begin{eqnarray}
\begin{array}{rlll}
\partial_t u &=\,\nabla\cdot \left( \tD \, \nabla u \right),
&\,x\in\Omega &\, t\in (0,T],\\[1mm]
u(x,t) &=\, 0,&\,x\in\partial\Omega,&\, t\in (0,T],\\[1mm]
u(x,0) &=\, u_0(x),&\,x\in\Omega,&
\end{array}
\end{eqnarray}
with domain $\Omega\subset\R^d\, (d\ge 1)$ and a symmetric positive definite matrix
$\tD (x)\in\R^{d\times d}$. We have the following stability results:
\begin{eqnarray}
&& \|u(t)\|_{L^2(\Omega)}\le \|u_0\|_{L^2(\Omega)},\quad t\in [0,T], \\[1mm]
&& u(t) \rightarrow 0 \mbox{ for } t\rightarrow\infty .
\end{eqnarray}
A Method of Lines approach (let's take finite differences for simplicity) yields
the system of ordinary differential equations
\begin{eqnarray}
\begin{array}{rll}
\partial_t U(t) &=\,A\,U(t),&\,t\in (0,T],\\[1mm]
U(0) &=\, U^0,&
\end{array}
\end{eqnarray}
where the vector $U(t)$ collects approximations at certain
spatial points. The matrix $A$ is symmetric negative definite
and therefore exhibits negative real eigenvalues - the values
$-k_i$ in the exponentials mentioned by Rosenbrock in
Fig.~\ref{fig:rosenbrock-paper-pic1}. An explicit Runge-Kutta
methods computes approximations $U_n\approx U(t_n)$ with
$t_n=nh,\, n\ge 1$ through
\begin{eqnarray}
\begin{array}{rll}
U_{n+1} &=\,R_{ERK} (h A) \,U_n,&\quad n=0,1,\ldots\\[2mm]
U_0 &=\, U^0,&
\end{array}
\end{eqnarray}
where
\begin{equation}
\label{eq:explicit-rk}
R_{ERK}(z) = 1+z+\ldots+\frac{z^p}{p!}+\sum_{i=p+1}^{s}\alpha_iz^i
= e^z + \Oh(z^{p+1}).
\end{equation}
The stability requirements for the semi-discretized solution
\begin{equation}
\|U_{n+1}\|_2\le \|U_n\|_2,\quad
U_n \rightarrow 0 \mbox{ for } n\rightarrow\infty,
\end{equation}
request $|R_{ERK}(z)|<1$ for $z$ along the negative real axis.
Due to the nature of the approximation (\ref{eq:explicit-rk}), small
time steps $h$ are necessary to guarantee stability. Moreover, the
finer the spatial discretization, the smaller the time steps must be,
showing that explicit methods, in general, are inefficient for the
solution of such kind of (stiff) problems.
\begin{figure}[ht]
\centering
\fbox{\includegraphics[width=0.85\textwidth]{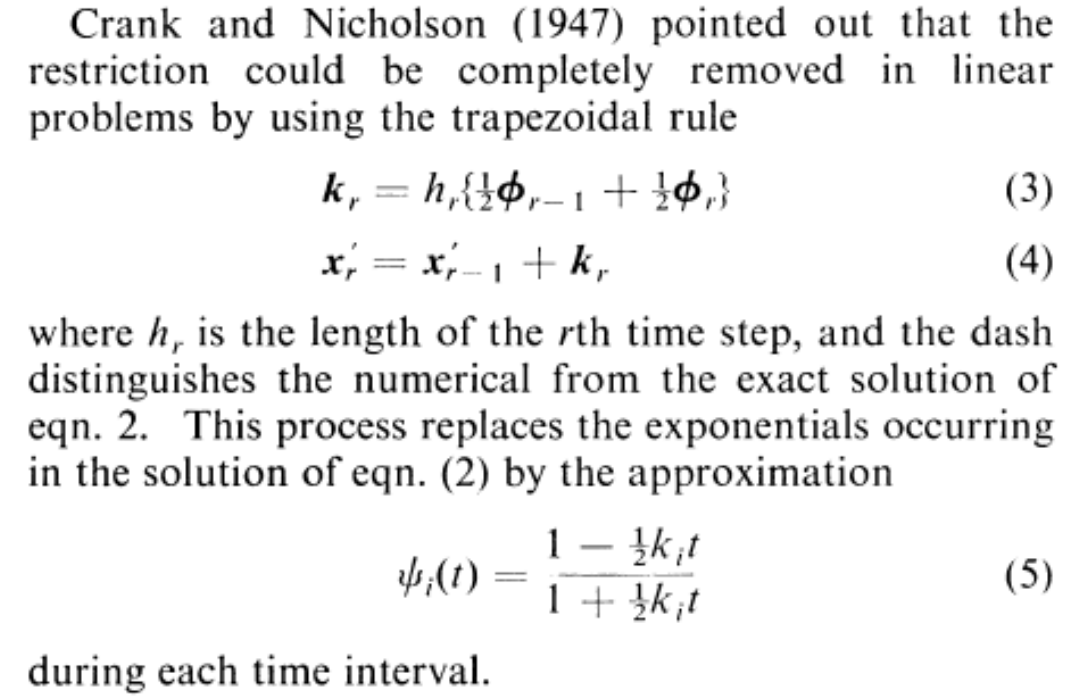}}
\parbox{13cm}{
\caption{\textcopyright\hspace{1mm}The Computer Journal, part of page 329 of \cite{Rosenbrock1963}.}
\label{fig:rosenbrock-paper-pic2}
}
\end{figure}

An alternative is the implicit Crank-Nicolson method
already proposed in 1947 \cite{CrankNicolson1947}, see
also Fig.~\ref{fig:rosenbrock-paper-pic2}. It reads
\begin{eqnarray}
\begin{array}{rll}
U_{n+1} &=\,R_{CN} (h A) \,U_n,&\quad n=0,1,\ldots\\[1mm]
U_0 &=\, U^0,&
\end{array}
\end{eqnarray}
with
\begin{equation}
R_{CN}(z) = \frac{1+z/2}{1-z/2} = e^z + \Oh(z^3).
\end{equation}
The method is unconditionally stable, since $|R_{CN}(z)|\le 1$
for all $z$ lying in the left complex half plane. However, the
damping properties at infinity are unsatisfactory. This lack has
been also mentioned by Rosenbrock: {\it The procedure given in eqns.
(3) and (4) has been widely used. It is perhaps not widely known,
however, that instability can arise even with this process when
the $\phi$ are non-linear functions of $x$. This is hardly surprising,
since $\psi_i(t)\rightarrow -1$ as $t\rightarrow\infty$, so that even in
a linear problem stability is only just maintained for large $t$.}

However, his main observation was that {\it when the functions $\phi$
are non-linear, implicit equations such as eqn. (3) can in general be
solved only by iteration. This is a severe drawback, as it adds
to the problem of stability, that of convergence of the iterative process.}
As consequence, he set up a generalized implicit process with
linear equations that can be solved rapidly and easily. How is it done?

\begin{figure}[ht]
\centering
\fbox{\includegraphics[width=0.85\textwidth]{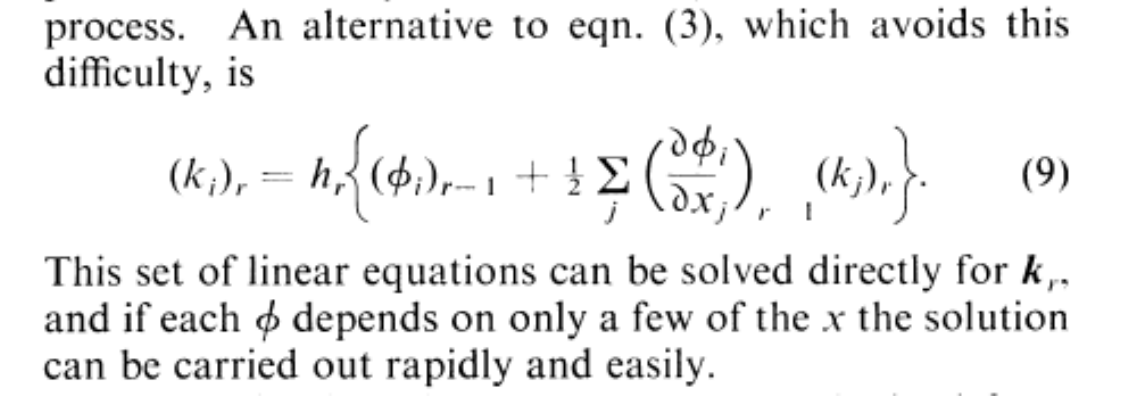}}
\parbox{13cm}{
\caption{\textcopyright\hspace{1mm}The Computer Journal, part of page 329 of \cite{Rosenbrock1963}.}
\label{fig:rosenbrock-paper-pic3}
}
\end{figure}
Let us now consider the autonomous system of nonlinear ordinary differential equations
\begin{eqnarray}
\label{eq:ode}
\begin{array}{rll}
\partial_t U(t) &=\,F(U(t)),&\,t\in (0,T],\\[1mm]
U(0) &=\, U^0.&
\end{array}
\end{eqnarray}
and apply the Crank-Nicolson method to it, resulting in
\begin{eqnarray}
\begin{array}{rll}
U_{n+1} &=\,U_n + \frac{h}{2} \left( F(U_n) + F(U_{n+1}) \right),&\quad n=0,1,\ldots\\[2mm]
U_0 &=\, U^0.
\end{array}
\end{eqnarray}
With $U_n$ as starting values, Newton's method to
approximate $U_{n+1}$ gives the sequence of linear equations
\begin{eqnarray}
\begin{array}{rll}
U_{n+1}^{(0)} &=\,U_n,\\[2mm]
\left( I - \frac{h}{2}F'(U_{n+1}^{(k)})\right) K_{n+1}^{(k+1)} &=\,
- \left( U_{n+1}^{(k)} - U_n - \frac{h}{2} \left( F(U_{n+1}^{(k)})
+ F(U_n)\right) \right),\\[4mm]
U_{n+1}^{(k+1)} &=\, U_{n+1}^{(k)} + K_{n+1}^{(k+1)},\quad k=0,1,\ldots .
\end{array}
\end{eqnarray}
The fundamental idea of Rosenbrock was to use {\it only one step} of
Newton's method, which reads for $k=0$
\begin{eqnarray}
\left( I - \frac{h}{2}F'(U_{n}) \right) K_{n+1}
=\, h \,F(U_{n}),
\end{eqnarray}
or equivalently
\begin{eqnarray}
K_{n+1} &=\, h \left( F(U_{n}) + \frac{1}{2}F'(U_{n})K_{n+1} \right).
\end{eqnarray}
In his paper, Rosenbrock did not mention how he derived his
formula (9), see Fig.~\ref{fig:rosenbrock-paper-pic3}, and left it to
the reader as an exercise.

\begin{figure}[ht]
\centering
\fbox{\includegraphics[width=0.85\textwidth]{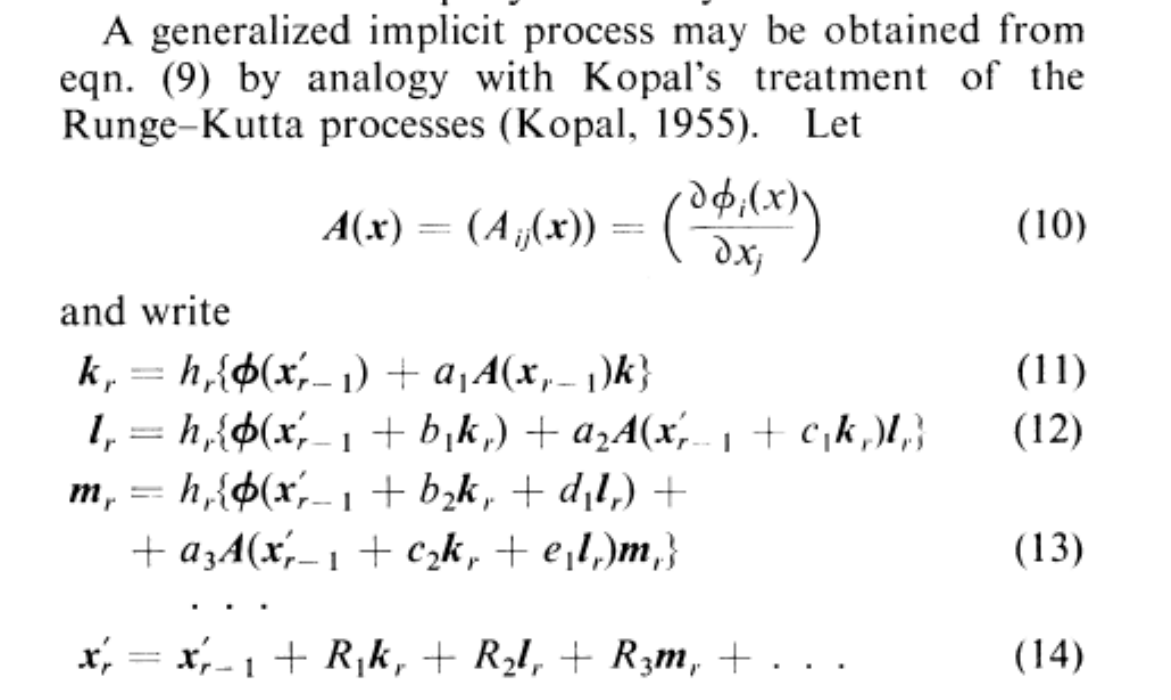}}
\parbox{13cm}{
\caption{\textcopyright\hspace{1mm}The Computer Journal, part of page 329 of \cite{Rosenbrock1963}.}
\label{fig:rosenbrock-paper-pic4}
}
\end{figure}

In a next step, he proposed to use Kopal's treatment of the
Runge-Kutta processes \cite{Kopal1955} to design a generalized implicit process
shown in Fig.~\ref{fig:rosenbrock-paper-pic4}.
Note that the Jacobian is evaluated at different solutions. Rosenbrock did
a consistency analysis: {\it By a straightforward but tedious calculation it
is possible to expand $x_r'-x_{r-1}'$ in eqn. (14) as a power series in $h_r$,
and to compare this with the Taylor's series.} He derived order conditions
for two stages up to order four. Finally, I summarize his findings:
\begin{itemize}
\item[\labelitemiv] There is no 2-stage third-order method with $R(\infty)=0$.
\item[\labelitemiv] He constructed a 2-stage third-order method with $R(\infty)=-0.8$.
\item[\labelitemiv] He constructed a 2-stage second-order method with $R(\infty)=0$.
\end{itemize}
Compared to the second-order Crank-Nicolson method, Rosenbrock found a 2-stage
second-order method with optimal damping property at infinity and only two linear equations that have to be solved in each time step. General $s$-stage Rosenbrock
or Rosenbrock-Runge-Kutta methods can be written in the (modern) form
\begin{eqnarray}
\begin{array}{rll}
\label{eq:rosenbrock-method}
\left( I-h\gamma_{ii}F'(U_n+\sum_{j=1}^{i-1} \delta_{ij}K_j)\right) K_i &=\,
h F(U_n+\sum_{j=1}^{i-1}\alpha_{ij}K_j),\quad i=1,\ldots,s,\\[3mm]
U_{n+1} &=\, U_n + \sum_{i=1}^s b_i K_i .
\end{array}
\end{eqnarray}
This formulation is the starting point for further improvements.

\section{The Improvement by Wanner}
Around 1973, Gerhard Wanner became interested in Rosenbrock schemes
and added his famous sum, $hF'(U_n)\sum_{j=1,\ldots,i-1}\gamma_{ij}K_j$, on
the right hand side in (\ref{eq:rosenbrock-method}),
keeping at the same time the Jacobian fixed, i.e., using $F'(U_n)$ for all
stages \cite[1977]{Wanner1977}. Rosenbrock-Wanner methods (short
ROW methods) with $s$ stages have the general form
\begin{eqnarray}
\label{eq:row-method}
\begin{array}{rll}
\left( I-h\gamma_{ii}F'(U_n)\right) K_i &=\,
h F(U_n+\sum_{j=1}^{i-1}\alpha_{ij}K_j) +
hF'(U_n)\sum_{j=1}^{i-1}\gamma_{ij}K_j,\\[3mm]
&\quad i=1,\ldots,s,\\[3mm]
U_{n+1} &=\, U_n + \sum_{i=1}^s b_i K_i .
\end{array}
\end{eqnarray}
In the spirit of Rosenbrock, they can be derived from
diagonally implicit Runge-Kutta methods (short DIRK methods),
applying only one
simplified Newton step with the Jacobian $F'(U_n)$ and using
already calculated stage values as starting values in the
calculation of subsequent stages. Applied
to (\ref{eq:ode}), the nonlinear system for the stage values $K_i$
of a DIRK method with lower triangular coefficient matrix $D=(d_{ij})$
reads
\begin{eqnarray}
\begin{array}{rll}
\label{eq:dirk-method}
K_i &=\, hF(U_n+\sum_{j=1}^{i}d_{ij}K_j),\quad i=1,\ldots,s.
\end{array}
\end{eqnarray}
One-step of a Newton-like iteration
\begin{eqnarray}
\begin{array}{rll}
&(I-hd_{ii}F'(U_n))\left( K_i-K_i^{(0)}\right) =
hF\left( U_n + \sum_{j=1}^{i-1} d_{ij}K_j+d_{ii}
K_j^{(0)}\right) - K_i^{(0)}
\end{array}
\end{eqnarray}
with starting values
\begin{eqnarray}
\begin{array}{rll}
&\,K_1^{(0)}=0,\quad
K_i^{(0)}=-\sum_{j=1}^{i-1}\frac{\gamma_{ij}}{d_{ii}}K_j,
\quad i=2,\ldots,s,
\end{array}
\end{eqnarray}
yields the ROW method (\ref{eq:row-method}) with
$\alpha_{ij}=d_{ij}-\gamma_{ij}$ and $\gamma_{ii}=d_{ii}$.
Compared to Rosenbrock's original form (\ref{eq:rosenbrock-method}),
the coefficients $\delta_{ij}$ were removed to avoid recalculations of Jacobians
and new coefficients $\gamma_{ij}$ were added to have enough parameters
for consistency and good stability properties.

A usual simplification is to set $\gamma_{ii}=\gamma$ for
all $i=1,\ldots,s$. In case of direct solvers,
it allows to reuse an LU-decomposition of the linear system
matrix $I-h\gamma F'(U_n)$. It also simplifies iterative solvers,
when matrix decompositions as preconditioners are used. To avoid
the matrix-vector multiplication,
one introduces $S_i=\sum_{j=1,\ldots,i}\gamma_{ij}K_j$
and solves
\begin{eqnarray}
\begin{array}{rll}
\label{eq:rosenbrock-method-tr}
\left( \frac{I}{h\gamma}-F'(U_n)\right) S_i &=\,
F(U_n+\sum_{j=1}^{i-1}a_{ij}S_j) + \sum_{j=1}^{i-1}\frac{c_{ij}}{h}S_j,\\[3mm]
&\quad i=1,\ldots,s,\\[3mm]
U_{n+1} &=\, U_n + \sum_{i=1}^s m_i S_i .
\end{array}
\end{eqnarray}
Defining the matrix $\Gamma\!=\!(\gamma_{ij})_{i,j=1}^s$ with
$\gamma_{ii}\neq 0$ for all $i$, the new
parameters are derived from
\begin{eqnarray}
\begin{array}{rll}
(a_{ij})_{i,j=1}^s &\!\!\!\!=\!\!\!\!& (\alpha_{ij})_{i,j=1}^s\Gamma^{-1},\quad
(c_{ij})_{i,j=1}^s\!=\text{diag} (\gamma_{11}^{-1},\ldots,\gamma_{ss}^{-1})-\Gamma^{-1},\\[2mm]
&& (m_1,\ldots,m_s) = (b_1,\ldots,b_s)\Gamma^{-1}.
\end{array}
\end{eqnarray}
Further generalizations to non-autonomous systems and systems of the
special multiplicative form
$M(t,U)\partial_tU=F(t,U)$, where $M$ might be singular, are also
possible \cite{HairerWanner1991,LubichRoche1990}.

So far, the Jacobian has to be computed at every time step, which
can be quite costly. Steihaug and Wolfbrandt
\cite[1979]{SteihaugWolfbrandt1979}
developed so-called W-methods that avoid exact Jacobians, i.e.,
$F'(U_n)\approx T_n$ with arbitrary matrix $T_n$. The idea is to keep
the Jacobian unchanged over several time steps while still ensuring stability.
Less restrictive
time lagged approximations of the form $T_n\approx F'(U_n)+\Oh(h)$ were
proposed by Scholz and Verwer \cite[1983]{ScholzVerwer1983}, see also
Scholz \cite[1978/79]{Scholz1978,Scholz1979}, and Kaps and Ostermann
\cite[1988/89]{Ostermann1988,KapsOstermann1989}.
Rahunanthan and Stanescu recently discussed high-order W-methods
\cite[2010]{RahunanthanStanescu2010}. They have been also applied
to optimal control problems in Lang and Verwer \cite[2013]{LangVerwer2013}.

The linear equations in (\ref{eq:rosenbrock-method-tr}) can be successively
solved. Order conditions were derived by applying the theory of Butcher series.
They can be found in Wolfbrandt \cite[1977]{Wolfbrandt1977},
Kaps \cite[1977]{Kaps1977}, N\o rsett and Wolfbrandt \cite[1979]{NorsettWolfbrandt1979}, and Kaps and Wanner \cite[1981]{KapsWanner1981}. Further details
and many more information are given in the books of Van der Houven \cite[1976]{Houwen1976} and Hairer and Wanner \cite[1991]{HairerWanner1991}.

For later use, I briefly recall the definition of a few fundamental stability
concepts. Applied to the famous scalar Dahlquist's test equation $y'=\lambda y$, $y_0=1$ with $\lambda\in\C$, a ROW method (as any other Runge-Kutta method) gives
$U_{n+1}=R(z)U_n$, where $z=\lambda h$. The function $R(z)$ is called the
stability function of the method and the set $S=\{ z\in\C:\, |R(z)|\le 1\}$ defines its stability domain. The exact solution of the test equation is stable in the entire
negative complex half plane $C^-=\{z:\,Re(z)\le 0\}$, and it seems likely that a numerical method should preserve this stability property. Dahlquist (1963) called a method \textit{A-stable} if $C^-\subset S$. If in addition $\lim_{z\rightarrow -\infty}R(z)=0$, the method is called \textit{L-stable} - a property that was introduced by Ehle (1969) and guarantees a fast damping for those $z$ having very large negative real parts. A convenient way to ensure L-stability for ROW methods is to require $\alpha_{si}+\gamma_{si}=b_i$ for $i=1,\ldots,s$, and $\sum_{j}\alpha_{sj}=1$. Such methods are called \textit{stiffly accurate}.
A weaker concept was established by Widlund (1967) who called a method
\textit{A($\alpha$)-stable} if the sector $S_\alpha=\{z:\,|\arg(-z)|\le\alpha,\,z\ne 0\}$ is contained in its stability region.

There are A-stable and L-stable ROW methods available. ROW methods share their
linear stability properties with (singly) diagonally implicit Runge-Kutta methods
introduced by Alexander \cite[1977]{Alexander1977}. The role of the stability
parameter $\gamma$ was studied in Wanner \cite[1980]{Wanner1980}. Continuous extensions
of Rosenbrock-type methods for a frequent graphical output were introduced by
Ostermann \cite[1990]{Ostermann1990}.

\section{Development of Rosenbrock-Wanner Methods}
{\bf First Solvers.} The theoretical investigation of Rosenbrock-Wanner
methods at the end of the 70s
laid the starting point for a broad and fast development of efficient solvers.
The fourth-order codes GRK4A and GRK4T proposed by Kaps and Rentrop \cite[1979]{KapsRentrop1979} were equipped with a step size control based
on embedded formulas of order three. The first one is A-stable whereas the
second is only A($89.3^o$)-stable, but comes with smaller truncation errors.
They were successfully tested on the 25 stiff test problems of
Enright, Hull and Lindberg \cite[1975]{EnrightHullLindberg1975}.
Gottwald and Wanner presented their back-stepping algorithm to improve
the reliability of Rosenbrock methods \cite[1981]{GottwaldWanner1981}.
Time-lagged Jacobian matrices and a modified Richardson extrapolation
for variable steps size control within a fourth-order A-stable
Rosenbock-Wanner scheme (named RKRMC) were tested by Verwer, Scholz, Blom, and
Louter-Nool \cite[1983]{VerwerScholzBlomLouterNool1983}. Further analysis and
experiments have been made by Verwer \cite[1982]{Verwer1982a,Verwer1982b}.
Implementation issues were discussed by Shampine \cite[1982]{Shampine1982}.
Veldhuizen investigated the D-stability of the Kaps-Rentrop methods \cite[1984]{Veldhuizen1984}. There are two options to estimate local errors:
embedding and Richardson extrapolation. Kaps, Poon, and Bui did a careful
comparison of these two strategies in \cite[1985]{KapsPoonBui1985}. The
performance of Rosenbrock methods for large scale combustion problems
discretized by the Method of Lines was investigated by Ostermann, Kaps,
and Bui \cite[1986]{OstermannKapsBui1986}.\\

\noindent {\bf Partitioned Methods.} It is often useful to split the solution
vector $U(t)$ into stiff and non-stiff components, say $U_s(t)$ and $U_n(t)$.
After an appropriate reordering of the original equations, this gives a partitioned
system
\begin{eqnarray}
\begin{array}{rll}
U'_s(t) &=\, F_s(U_s(t),U_n(t)),\quad U_s(0)=U^{0}_s,\\[1mm]
U'_n(t) &=\, F_n(U_s(t),U_n(t)),\quad U_n(0)=U^{0}_n.
\end{array}
\end{eqnarray}
Now it is quite natural to apply a Rosenbrock-type scheme to the
stiff part and an explicit Runge-Kutta method to the non-stiff part.
Rentrop combined an A-stable Rosenbrock (3)4-pair with a common
(4)5-Runge-Kutta-pair and studied strategies for stiffness detection
in \cite[1985]{Rentrop1985}. A drawback of such an approach is the
occurrence of additional coupling conditions which usually does not
allow the simple combination of two favourite schemes. An alternative is
to use the setting of W-methods to directly incorporate the partitioning
on the level of the Jacobian calculation, e.g., only take into
account derivatives of $F_s$ and drop the other ones. Such methods were
analysed by Strehmel, Weiner, and Dannehl
\cite[1990]{StrehmelWeinerDannehl1990} under the heading
partitioned linearly implicit Runge-Kutta methods including ROW- and
W-methods. Later on, Wensch designed an eight-stage fourth-order
partitioned Rosenbrock method for multibody systems in index-3 formulation
\cite[1998]{Wensch1998}.

The partitioning can be also used to set up multirate schemes, where
different step sizes for active and latent components are explicitly
introduced in the discretization. In G\"unther and Rentrop
\cite[1993]{GuentherRentrop1993}, multirate Rosenbrock-Wanner methods
were used for the simulation of electrical networks. One general shortcoming
of multirate methods is the coupling between the components by interpolating and extrapolating state variables. Stability of multirate Rosenbrock methods
were studied in Savcenco \cite[2008/09]{Savcenco2008,Savcenco2009}
and Kuhn and Lang \cite[2014]{KuhnLang2014}.\\

\noindent {\bf Differential-Algebraic Equations.} In the late 80s,
Rosenbrock methods were also applied to differential-algebraic
equations (DAEs) of index one:
\begin{eqnarray}
\begin{array}{rll}
\label{eq:daes}
U'(t) &=\, F(U(t),Z(t)),\quad U(0)=U^0,\\[1mm]
0 &=\, G(U(t),Z(t)),\quad Z(0)=Z^0,
\end{array}
\end{eqnarray}
where it is assumed that $(\partial_Z G)^{-1}$ exists and is bounded
in a neighbourhood of the solution. The main idea used by Roche
\cite[1988]{Roche1988} is to add $\varepsilon Z'(t)$ on the left hand
side of the second equation and consider the DAE (\ref{eq:daes}) as a
limit case of the stiff singular perturbation problem for
$\varepsilon\rightarrow 0$. This limit typically destroys the
classical order of the Rosenbrock methods and gives rise to a new
consistency theory derived by means of a modified Butcher-like tree
model for the $U$- and $Z$-components. Note that the Kaps-Rentrop
methods from \cite{KapsRentrop1979} drop down to order two when applied
to (\ref{eq:daes}). Similar observations have been made earlier by Verwer
\cite[1982]{Verwer1982b}. Two new ROW-methods (named DAE34 and RKF4DA)
with stepsize control and
an index-1 monitor were proposed and tested by Rentrop, Roche and
Steinebach \cite[1989]{RentropRocheSteinebach1989}.

A desirable property when solving stiff or differential-algebraic equations
is to have an L-stable method, i.e., a method with $R(\infty)=0$. This is always
the case for stiffly accurate Rosenbrock methods which approximate the
algebraic component $Z$ of the extreme DAEs, $U'=1$ and $0=G(U,Z)$, through one
simplified Newton iteration. This nicely meets the original idea of Rosenbrock.
In their book, Hairer and Wanner \cite[1991]{HairerWanner1991} constructed the famous
stiffly accurate six-stage fourth-order Rosenbrock solver {\sc Rodas} with
an embedded method of order three. Special index-2 DAEs were treated in
Lubich and Roche \cite[1990]{LubichRoche1990} and results for index-3 multibody
systems can be found in Wensch \cite[1998]{Wensch1998}. G\"unther, Hoschek, and
Rentrop constructed special index-2 Rosenbrock methods for
electric circuit simulations \cite[2000]{GuenterHoschekRentrop2000}.
Recently, Jax and Steinebach \cite[2017]{JaxSteinebach2017} introduced a new type of ROW methods for solving DAEs of the form (\ref{eq:daes}). Taking ideas from W-methods, they allow arbitrary approximations to Jacobian entries resulting from the differential part.\\

\noindent {\bf Extrapolation.} An interesting, general approach to construct higher order methods for
differential as well as differential-algebraic equations is to use
extrapolation. Deuflhard and Nowak \cite[1987]{DeuflhardNowak1987} proposed
to extrapolate the linearly implicit Euler discretization (as the simplest
Rosenbrock method) to solve chemical reaction kinetics and electric circuits
and implemented the well-known variable-order {\sc Limex} code with step size
control. They also provided the impetus for Lubich to explain the error
behaviour of such methods by perturbed asymptotic analysis
\cite[1989]{Lubich1989}.\\

\noindent {\bf B-Convergence and Order Reduction.} One-step methods and so
Rosenbrock schemes suffer from order reduction, especially when they are applied
to nonlinear parabolic partial differential equations. Sharp error estimates
showing fractional orders of convergence for Rosenbrock and
W-methods were first established by Lubich, Ostermann, and Roche \cite[1993/95]{OstermannRoche1993,LubichOstermann1995}. This phenomenon
is related to the B-convergence of linearly implicit methods studied
by Strehmel and Weiner \cite[1987]{StrehmelWeiner1987}. Barriers for
the order of B-convergence were given by Scholz \cite[1989]{Scholz1989}.
In their book, Strehmel and Weiner \cite[1992]{StrehmelWeiner1992} gave
convergence results for spatial discretizations of semilinear parabolic
equations with constant operator and a Lipschitz continuous non-linearity.
However, the B-convergence technique does not give the sharp fractional
temporal convergence rates. It is now much better understood than
before why (lower) fractional orders occur. This reduction is not induced
by lack of smoothness of the solution but rather by the presence of powers
of the spatial differential operators in the local truncation error.
Concerning W-methods, the order reduction is more severe compared
with Rosenbrock methods. Loss of accuracy happens long before stability
is affected. Fortunately, there are additional consistency conditions
that imply also higher order of convergence as shown in
Lubich and Ostermann \cite[1995]{LubichOstermann1995}.

Using this theoretical framework, new methods
were constructed. Steinebach improved the {\sc Rodas} code and designed
his stiffly accurate {\sc Rodasp} scheme, which satisfies the new
conditions for linear parabolic problems to reach order four. It was
successfully applied  to forecast transport in rivers, see Steinebach
and Rentrop \cite[2001]{SteinebachRentrop2001}. New order-three methods
with three, {\sc Ros3p}, and four stages, {\sc Ros3pl}, were constructed
in Lang and Verwer \cite[2001]{LangVerwer2001} and Lang and Teleaga
\cite[2008]{LangTeleaga2008}, respectively. The latter one is stiffly
accurate and therefore suitable for differential-algebraic equations.
It also satisfies the condition of a W-method with $\Oh(h)$-disturbance
of the Jacobian, which makes numerical differentiation for its entries
less sensitive with respect to roundoff errors. A bunch of newly
designed third-order Rosenbrock W-methods for partial
differential-algebraic equations was published in Rang and Angermann
\cite[2005]{RangAngermann2005}. Further improved ROW methods can be
found in Rang \cite[2014/15]{Rang2014,Rang2015}.\\

\noindent {\bf Exponential Rosenbrock-type Methods.} Exponential integrators
are based on a continuous linearization of the nonlinearity $F(U(t))$ along
the numerical solution. This gives the linearized system
\begin{equation}
U'(t) = F'(U_n)U(t) + G_n(U(t)),\quad
G_n(U(t))=F(U(t))-F'(U_n)U(t).
\end{equation}
Exponential Rosenbrock methods make direct use of $J_n:=F'(U_n)$ and $G_n(U(t))$.
Hochbruck, Ostermann, and Schweitzer \cite[2009]{HochbruckOstermannSchweitzer2009}
considered the following class of methods (here for variable time steps $h_n$):
\begin{equation}
\begin{array}{rll}
U_{ni} &=& e^{c_ih_nJ_n}\,U_n+h_n\sum_{j=1}^{i-1}a_{ij}(h_nJ_n)g_n(U_{nj})
,\quad i=1,\ldots,s,\\[3mm]
U_{n+1} &=& e^{c_ih_nJ_n}\,U_n+h_n\sum_{i=1}^{s}b_{i}(h_nJ_n)g_n(U_{ni}).
\end{array}
\end{equation}
A key point is the efficient approximation of the matrix exponential times
a vector by Krylov subspace methods or methods based on direct polynomial
interpolation. An interpolation method with real Leja points was tested
by Caliari and Ostermann \cite[2009]{CaliariOstermann2009} and showed a great
potential for problems with large advection in combination with moderate diffusion
and mildly stiff reactions. Higher order and parallel exponential Rosenbrock
methods were proposed by Luan and Ostermann \cite[2014/16]{LuanOstermann2014,LuanOstermann2016}.\\

\noindent {\bf Miscellaneous.} Rosenbrock methods offer a simple usage due
to their linear structure. Methods up to order four perform well for low
and medium tolerances and work competitive in many applications. The code
ode23s in the {\sc Matlab Ode Suite} is a typical Rosenbrock scheme, see
Shampine and Reichelt \cite[1997]{ShampineReichelt1997}. The Krylov-W-code
{\sc Rowmap} based on the Rosenbrock method {\sc Ros4} of Hairer and Wanner
has demonstrated its efficiency for large stiff systems. Numerical tests
were performed in Weiner, Schmitt, and Podhaisky
\cite[1997]{WeinerSchmittPodhaisky1997}. Rosenbrock methods are the
numerical kernel in the adaptive multilevel PDAE-solver {\sc Kardos}, which
is a well running working horse for a broad range of real-life applications,
see Lang \cite[2000]{Lang2000}. Combined with a linearized error transport
equation based on first variational principles, they can be accompanied with
a cheap global error estimation and control through tolerance proportionality.
Such strategies were investigated in Lang and Verwer
\cite[2007]{LangVerwer2007} for initial value problems
and in Debrabant and Lang \cite[2015]{DebrabantLang2015} for
semilinear parabolic equations.
Last but not least, a Rosenbrock code is listed
in the second edition of {\it Numerical Recipes} by Press, Teukolsky,
Vetterling, and Flannery \cite[1996]{PressTeukolskyVetterlingFlannery1996}.

A lot of basic information about Rosenbrock methods can be found in the
books by Hairer and Wanner \cite[1991]{HairerWanner1991} and Strehmel and Weiner \cite[1992]{StrehmelWeiner1992}. Newer developments are highlighted in
Strehmel, Weiner, and Podhaisky \cite[2012]{StrehmelWeinerPodhaisky2012}.
A tremendous source of further interesting material are the proceedings
of the numerous NUMDIFF-conferences held at the Martin Luther University
Halle-Wittenberg since the early 1980s.

\section{Two-Step Rosenbrock-Peer and W-Methods}
As explained above, Rosenbrock methods may suffer from order reduction for
very stiff problems. A closer inspection reveals that the low stage order
(the first stage value is computed by the linearly implicit Euler scheme)
is one of the reasons. To raise the stage order substantially, Podhaisky,
Schmitt, and Weiner \cite[2002]{PodhaiskyWeinerSchmitt2002,PodhaiskySchmittWeiner2002} studied a new class of linearly implicit two-step methods,
where the previously computed stage values are taken into account. Such
$s$-stage two-step W-methods have the form
\begin{equation}
\begin{array}{rll}
Y_{ni} &=& U_n + h_n\sum_{j=1}^{s}a_{ij}U_{n-1,j}+
h_n\sum_{j=1}^{i-1}\tilde{a}_{ij}U_{nj},\\[2mm]
(I-\gamma h_nT_n) U_{ni} &=& F(Y_{ni})+h_nT_n\sum_{j=1}^{s}\gamma_{ij}U_{n-1,j}+
h_nT_n\sum_{j=1}^{i-1}\tilde{\gamma}_{ij}U_{nj},\\[3mm]
&& i=1,\ldots,s,\\[2mm]
U_{n+1} &=& U_n + h_n\sum_{i=1}^{s} \left( b_iU_{ni}+v_iU_{n-1,i}\right).
\end{array}
\end{equation}
Observe that $a_{ij}\!=\!\gamma_{ij}\!=\!v_i\!=\!0$ recovers classical
one-step ROW and W-methods. The special setting $\tilde{a}_{ij}\!=\!\tilde{\gamma}_{ij}\!=\!0$ treated in \cite{PodhaiskySchmittWeiner2002} allows to compute the stage values
$U_{ni}$ in parallel. Higher order parallel methods were studied
by Jackiewicz, Podhaisky, and Weiner \cite[2004]{JackiewiczPodhaiskyWeiner2004}.
Computer architectures of workgroup servers having shared memory for
quite a few processors are particularly suitable for these methods which
have been designed for the solution of large stiff systems in combination
with Krylov techniques. Methods with favorable stability properties have been
constructed with stage order $q\!=\!s$ and order $p\!=\!s$ for $s\!\le\!4$.
All methods are competitive with state-of-the-art codes for stiff ODEs.

Within the class of two-step methods, Podhaisky, Weiner,
and Schmitt \cite[2005]{PodhaiskyWeinerSchmitt2005} also constructed $s$-stage
methods where all stage values have the stage order $q\!=\!s\!-\!1$. They
considered the following methods:
\begin{equation}
\begin{array}{rll}
(I-\gamma h_nT_n) U_{ni} &=& \sum_{j=1}^{s}b_{ij}U_{n-1,j}
+ h_n\,\sum_{j=1}^{s} a_{ij} \left( F(U_{n-1,j}) - T_nU_{n-1,j}\right)\\[3mm]
&& + h_nT_n \sum_{j=1}^{i-1}g_{ij}U_{nj},\quad i=1,\ldots,s.
\end{array}
\end{equation}
Here, $U_{ns}\approx U(t_{n+1})$ and the matrix $T_n$ is supposed to be an
approximation to the Jacobian $F'(U(t_n))$ for stability reasons. The method
is treated as a W-method, i.e., the order conditions are derived for arbitrary
$T_n$. Due to their two-step and linear structure, the methods are called
two-step Rosenbrock-Peer methods, where {\it peer} refers to the fact that
all stage values have now one and same order. The methods constructed
in \cite{PodhaiskyWeinerSchmitt2005} for $s=4,\ldots,8$ are zero-stable for
arbitrary step size sequences and L($\alpha$)-stable with large $\alpha$.
For constant time steps, these methods have order $s$. Numerical experiments
showed no order reduction and an efficiency superior to the fourth-order
{\sc Rodas} for more stringent tolerances.

With this property, peer methods commend themselves as time-stepping schemes
for the solution of time-dependent partial differential equations. So they
have been implemented in the already mentioned finite element software
package {\sc Kardos}, see Gerisch, Lang, Podhaisky, and Weiner \cite[2009]{GerischLangPodhaiskyWeiner2009} and Schr\"oder, Gerisch,
and Lang \cite[2017]{SchroederGerischLang2017}. They also performed well for
compressible Euler equations, demonstrated in Jebens, Knoth, and
Weiner \cite[2012]{JebensKnothWeiner2012}, for shallow-water equations,
reported in Steinebach and Weiner \cite[2012]{SteinebachWeiner2012}, and
for more complex fluid dynamics problems, see Gottermeier and Lang
\cite[2009/10]{GottermeierLang2009,GottermeierLang2010}. More
recently, linearly implicit two-step Peer methods of Rosenbrock-type
have shown their reliability, robustness, and accuracy for large eddy
and direct numerical simulations for turbulent unsteady flows in
Massa, Noventa, Lorini, Bassi, and Ghidoni
\cite[2018]{MassaNoventaLoriniBassiGhidoni2018}.

\section{Summary}
The idea of Rosenbrock is still alive. Avoiding the (often cumbersome)
solution of nonlinear equations has not lost its attractiveness and
significance over the years. The successive solution of linear equations
is still a valuable option to efficiently solve systems of differential,
differential-algebraic or partial differential equations. Classical
one-step Rosenbrock-Wanner methods up to order four have demonstrated
their good performance for low and medium tolerances. The new class of
two-step Rosenbrock-Peer methods allows the construction of even higher
order methods that overcome the disadvantage of order reduction and still
exhibit good stability properties. Recent numerical experiments
with higher tolerances are very promising.

There is still an ongoing research activity in the field of Rosenbrock
methods. A recent search in the {\sc Scopus} data base gave $753$ documents.
One of the last entries is about {\it Strong Convergence Analysis of the Stochastic
Exponential Rosenbrock Scheme for the Finite Element Discretization of
Semilinear SPDEs Driven by Multiplicative and Additive Noise} by
Mukam and Tambue \cite[2018]{MukamTambue2018}. This brings me to my final
remark. In view of the numerous contributions to Rosenbrock schemes, I
would like to apologize in advance to those who have made significant further
contributions to the topic but were not mentioned in my overview. I am
prepared to receive your emails.

\section{Acknowlegdement}
I would like to thank R\"udiger Weiner for careful
reading of a first version of the manuscript. His suggestions very much
helped to improve the article. I also thank the reviewers for their useful
remarks.

\bibliographystyle{plain}
\bibliography{bibrow}

\end{document}